\documentstyle[12pt]{article}
\advance \textheight by \topskip
\makeatletter
\def\eqnarray{\stepcounter{equation}\let\@currentlabel=\theequation
\global\@eqnswtrue
\global\@eqcnt\z@\tabskip\@centering\let\\=\@eqncr
$$\halign to \displaywidth\bgroup\@eqnsel\hskip\@centering
  $\displaystyle\tabskip\z@{##}$&\global\@eqcnt\@ne
  \hfil$\displaystyle{\hbox{}##\hbox{}}$\hfil
  &\global\@eqcnt\tw@ $\displaystyle\tabskip\z@
  {##}$\hfil\tabskip\@centering&\llap{##}\tabskip\z@\cr}
\@addtoreset{equation}{section}
  \def\theequation{\thesection.\arabic{equation}}
\makeatother
\def\inprod{\mathop{\kern -0.05em\raise -0.1em\hbox{%
  \vrule height 0.03em width 0.6em depth 0em%
  \vrule height 0.7em width 0.03em depth 0em}\kern 0.1em}\nolimits}

\begin{document}

\begin{titlepage}
\hbox to \hsize{\hfil q-alg/9501017}
\hbox to \hsize{\hfil IHEP 94--145}
\hbox to \hsize{\hfil December, 1994}
\vfill
\large \bf
\begin{center}
Boundary values as Hamiltonian variables. II. \\
Graded structures
\end{center}
\vskip 1cm
\normalsize
\begin{center}
{\bf Vladimir O. Soloviev\footnote{e--mail:
vosoloviev@mx.ihep.su}}\\
{\small Institute for High Energy Physics,\\
 142284 Protvino, Moscow Region, Russia}
\end{center}
\vskip 2.cm
\begin{abstract}
\noindent
It is shown that new Poisson brackets proposed in Part I of this
work (hep-th/9305133) arise naturally in an
extension of the formal variational calculus incorporating divergences. The
linear spaces of local functionals, evolutionary vector fields, functional
forms, multi-vectors and differential operators become graded with respect to
divergences. Bilinear operations, such as action of vector fields onto
functionals, commutator of vector fields, interior product of  forms and
vectors and the Schouten-Nijenhuis bracket are compatible with the grading.
A definition of the adjoint graded operator is proposed and skew-adjoint
operators are constructed with the help of boundary terms.
Fulfilment of the Jacobi identity for the new Poisson brackets is shown to be
equivalent to the vanishing of the Schouten-Nijenhuis bracket for Poisson
bivector with itself. The simple procedure for testing this condition
proposed by P. Olver is applicable with a minimal modification. It is
demonstrated that the second structure of the Korteweg-de Vries equation is
not Hamiltonian with respect to the new brackets.
\end{abstract}
\vfill
\end{titlepage}

\def\d{\mbox{\sf d}}
\def\Rn{\rm R^n}

\section {Introduction}
In the first part of this work\cite{s1} (below simply I) field theory Poisson
brackets have been constructed which fulfil the Jacobi identity under
arbitrary boundary conditions.
The purpose of this paper is to show that these brackets replace
the standard ones in the extension of the formal variational calculus
onto total divergences. This extension consists in the introduction of a new
grading for linear spaces of local functionals,
vector fields, functional forms, multi-vectors and differential operators.
Nonstandard Rules 4.2 and 5.4 of I for the treatment of distributions
on closed
domains can be better understood in this context. The former
can now be treated as a definition of the interior product
of 1-forms and multi-vectors or, vice versa, 1-vectors and differential forms,
and the latter as a manifestation of
compatibility of the bilinear operations with the grading.

We show that the Jacobi identity for the new Poisson brackets can be
verified without lengthy calculations of I. Its fulfilment is
equivalent to the vanishing of the Schouten--Nijenhuis bracket of the
Poisson bivector with itself. And in its turn this condition can be easily
tested along with the procedure proposed by  Olver\cite{olv} with a
minimal modification of it. So, more attention than in I is paid here to
nonultralocal Hamiltonian operators with nonconstant coefficients. It occurs
that not all operators which are Hamiltonian with respect to the standard
brackets remain Hamiltonian in relation to the new brackets. For example, the
second structure of the Korteweg--de Vries equation is not Hamiltonian.

As a rule, we use the same notations as in I, except that we
change the notation for the Fr\'echet derivative  from $D_{f}$ to $f'$.
We still find convenient to represent integrals over finite domain $\Omega$
through integrals over the infinite space $\Rn$ by inserting
into all integrands  the characteristic function
$\theta_{\Omega}$. Then the formalism seems closer to the standard
formal variational calculus where local functionals and functional forms
are defined modulo divergences. But the formal divergences that we
discard here are integrated to zero under arbitrary
conditions on the boundary of the finite domain, whereas real
divergences are incorporated into graded structures.
All the operations introduced below are compatible both with
discarding formal total divergences
(if one object is a formal divergence than the result of operation
is also formal divergence) and with grading
(i.e., the same is valid for real divergences).
Extension of the space of differential operators by admitting their grading
permits
to use the concept of adjoint operator. So, skew-adjoint operators can now
be constructed and the Poisson bracket formulas become shorter,
though their content is the same. Nevertheless, in the proof of the Jacobi
identity we prefer to use the old notations  for easier comparison with
partial proofs of I.

The content of the paper is as follows.
Sections 2 and 3 are devoted to the introduction of the grading for the basic
structures of the formal variational calculus: local functionals,
vector fields and differential forms. Section 4 deals with graded
differential operators and their adjoints. In Section 5 we discuss
multi-vectors and the  Schouten-Nijenhuis bracket. Section 6
is devoted to concepts of the Hamiltonian formalism with the proof of
Jacobi identity postponed until Section 7. At last, in Section 8
we consider two examples: the second structure of Korteweg-de Vries
equation and the 2-dimensional flow of the ideal fluid. A short summary
is given in Conclusion.

\section{Local functionals and evolutionary vector fields}
Let us start with notions from the theory of graded spaces as they are
given in Ref.~\cite{dorf}.
A {\it  grading} in linear space $L$ is a decomposition of it into direct sum
of subspaces, with a special value of some function $p$ (grading function)
assigned to all the elements of any subspace.

Below the function $p$ takes its values in the set of all positive
multi-indices $J=(j_1,\dots,j_n)$ and so,
\[
L=\bigoplus\limits_{J=0}^{\infty} L^{\langle J\rangle}.
\]
Elements of each subspace are called homogeneous.

A bilinear operation $x,y\mapsto x\circ y$, defined on $L$, is said to be
{\it compatible with the grading} if the product of any homogeneous elements
is also homogeneous, and if
\[
p(x\circ y)=p(x)+p(y).
\]

Now  turn to concrete structures.
The space of local functionals $\cal A$ has already been defined in I.
Here we will call the
expression given in Definition 2.1 of I
the {\it canonical form of a local functional}. We formally extend
that definition by allowing local functionals to be written as follows
\[
F=\sum_{J=0}^{\infty}\int D_J\theta_{\Omega}f^{\langle J\rangle}
\biggl(\phi_A(x),
D_K\phi_A(x)\biggr)d^nx=\sum\int\theta^{(J)}f^{\langle J\rangle},
\]
where in accordance with the previous definition only a finite number
of terms is allowed. Here and below we simplify the notation
for derivatives of $\theta$ and remove $\Omega$.
Of course, any functional of such a form can be transformed
to the form used in I through integration by parts
\[
F=\int\theta_{\Omega}f=
\int\limits_{\Omega}f,
\]
where
\[
f=\sum(-1)^{|J|}D_Jf^{\langle J\rangle}.
\]
So, the formal integration by parts over infinite
space $\Rn$ evidently changes the grading.
It will be clear below that the general
situation is from one side compatibility of all bilinear operations
with the grading and from the other side with formal integration by parts.
So, basic objects (local functionals etc.) are defined as equivalence
classes modulo formal divergences (i.e., divergences of expressions
containing $\theta$-factors) and the unique
decomposition into homogeneous subspaces with fixed grading function
can be made only for representatives of these classes.

We call  expressions of the form
\[
\psi=\sum\int\theta^{(J)}
D_K\psi^{\langle J\rangle}_A\frac{\partial}{\partial\phi_A^{(K)}}
\]
the {\it evolutionary vector fields}.
Value of the evolutionary vector field on a local functional is given
by the expression
\[
\psi F=\sum\int\theta^{(I+J)}
D_K\psi^{\langle J\rangle}_A\frac{\partial f^{\langle I
\rangle}}{\partial\phi_A^{(K)}}.
\]
In principle, this formula can be understood  as a definition
but we can interpret it also as a consequence of
the standard relation
\[
\frac{\partial\phi_A(y)}{\partial\phi_B(x)}=
\delta (x,y)\delta_{AB}
\]
and Rule 5.4 of I.
It is a straightforward calculation to check that this operation is
compatible with the formal integration by parts, i.e.
\[
\psi{\rm Div}(f)={\rm Div}(\psi f),
\]
as it is in the standard formal variational calculus. This relation is,
of course, valid for integrands.

It is easy to see that the evolutionary vector field with coefficients
\[
\psi^{\langle J\rangle}_A=\sum
\biggl( D_L\xi_B^{\langle I\rangle} \frac{\partial\eta_A^{\langle J-I\rangle
}}{\partial \phi_B^{(L)}}-
D_L\eta_B^{\langle I\rangle} \frac{\partial\xi_A^{\langle J-I\rangle}}
{\partial \phi_B^{(L)}}\biggr)
\]
can be considered as the {\it  commutator of the evolutionary vector fields}
 $\xi$ and $\eta$
\[
\psi F=[\xi,\eta]F=\xi(\eta F)-\eta(\xi F),
\]
with the Jacobi identity fulfilled for the commutator operation,
and so these vector fields form a Lie algebra.

\section{Differentials and functional forms}
The {\it differential of a local functional}
is simply the first variation of it
\[
\d F=\sum\int\theta^{(J)}
\frac{\partial f^{\langle J\rangle}}{\partial\phi_A^{(K)}}\delta\phi_A^{(K)},
\]
where here and below $\delta\phi_A^{(K)}=D_K\delta\phi_A$.
It can also be expressed through
the Fr\'echet derivative (Definition 2.13 of I)
or through the higher Eulerian operators (Definition 2.4 of I)
\[
\d F=\sum\int\theta^{(J)}{f^{\langle J\rangle}}'(\delta\phi)=
\sum\int\theta^{(J)}
D_K\biggl( E^K_A(f^{\langle J\rangle})\delta\phi_A\biggr) .
\]
This differential is a special example of functional 1-form.
A general functional 1-form can be written as
\[
\alpha = \sum\int\theta^{(J)}\alpha ^{\langle J\rangle}_{AK}
\delta\phi_A^{(K)}.
\]
Of course, the coefficients are not unique since we can make
formal integration by parts. Let us call the following
expression the {\it canonical form of a functional 1-form}
\[
\alpha=\sum\int\theta^{(J)}\alpha^{\langle J\rangle}_A\delta\phi_A.
\]
Analogously, we can define {\it functional $m$-forms}
as integrals or equivalence classes modulo formal
divergences of vertical forms
\[
\alpha =\frac{1}{m!}
\sum\int\theta^{(J)}\alpha^{\langle J\rangle}_{A_1K_1,\dots,A_mK_m}\delta
\phi_{A_1}^{(K_1)}\wedge\dots\wedge\delta\phi_{A_m}^{(K_m)}.
\]
Define the {\it pairing} (or the {\it interior product})
of an evolutionary vector field and 1-form as
\begin{equation}
\alpha (\xi)=\xi\inprod\alpha=\sum
\int\theta^{(I+J)}\alpha ^{\langle J\rangle}_{AK}D_K\xi_A^{\langle I\rangle}.
\label{eq:pairing}
\end{equation}
The interior product of an evolutionary vector field and functional $m$-form
will be
\[
\xi\inprod\alpha=
\frac{1}{m!}\sum(-1)^{i+1}
\int\theta^{(J+I)}\alpha^{\langle J\rangle}_{A_1K_1,\dots,A_mK_m}
D_{K_i}\xi_{A_i}^{\langle I\rangle}\delta
\phi_{A_1}^{(K_1)}\wedge\dots
\]
\[
\dots\wedge\delta\phi_{A_{i-1}}^{(K_{i-1})}
\wedge\delta\phi_{A_{i+1}}^{(K_{i+1})}\wedge\dots\wedge
\delta\phi_{A_m}^{(K_m)}.
\]
Then the value of $m$-form on the $m$ evolutionary vector fields
will be defined by the formula
\[
\alpha (\xi_1,\dots,\xi_m)=\xi_m\inprod\dots\xi_1\inprod\alpha.
\]
It can be  checked by straightforward calculation that
\[
{\rm Div}(\alpha) (\xi_1,\dots,\xi_m)=
{\rm Div}(\alpha (\xi_1,\dots,\xi_m)).
\]
The {\it differential of the $m$-form} given as
\[
\d\alpha =\frac{1}{m!}
\sum
\int\theta^{(J)}\frac{\partial\alpha^{\langle J
\rangle}_{A_1K_1,\dots,A_mK_m}}
{\partial\phi_A^{(K)}}\delta\phi_A^{(K)}\wedge\delta
\phi_{A_1}^{(K_1)}\wedge\dots
\wedge\delta\phi_{A_m}^{(K_m)},
\]
satisfies standard properties
\[
{\d}^2=0
\]
and
\[
\d\alpha(\xi_1,\dots,\xi_{m+1})=
\sum\limits_i(-1)^{i+1}\xi_i\alpha(\xi_1,\dots,
\hat\xi_i,\dots,\xi_{m+1})+
\]
\[
+\sum\limits_{i<j}(-1)^{i+j}\alpha([\xi_i,\xi_j],\xi_1,\dots,\hat\xi_i,\dots,
\hat\xi_j,\dots,\xi_{m+1}).
\]
The {\it Lie derivative} of a functional form $\alpha$
along an evolutionary vector field $\xi$ can be introduced by the standard
formula
\[
L_{\xi}\alpha=\xi\inprod\d\alpha+\d \biggl(\xi\inprod\alpha\biggr).
\]

\section{Graded differential operators and their adjoints}
We call linear matrix differential operators of the form
\[
\hat I=\sum_{J\ge 0}\theta^{(J)}
\sum_{N=0}^{N_{max}}I^{\langle J\rangle N}_{AB}D_N
\]
the {\it graded differential operators}.

Let us call the linear differential operator $\hat I^{\ast}$ the {\it
adjoint} to $\hat I$ if for an
arbitrary set of smooth functions $f_A$, $g_A$
\[ \sum\limits_{A,B}\int f_A\hat I_{AB}g_B= \sum\limits_{A,B}\int g_A\hat
I^{\ast}_{AB}f_B.  \] For coefficients of the adjoint operator we can derive
the expression \begin{equation} I^{\ast\langle J\rangle
M}_{AB}=\sum\limits_{K=0}^{K_{max}} \sum\limits_{L=0}^{min(K,J)}
(-1)^{|K|}{K\choose L}{K-L\choose M}
D_{K-L-M}I^{\langle J-L\rangle K}_{BA}.\label{eq:adj}
\end{equation}
It is easy to check that the relation
\[
\hat I(x)\delta(x,y)=\hat I^{\ast}(y)\delta(x,y)
\]
follows from Rule 4.2 of I. For example, we have
\begin{equation}
\biggl(\frac{\partial}{\partial x^i}+\frac{\partial}{\partial y^i}\biggr)
\delta (x,y)=-\theta^{(i)}\delta (x,y).\label{eq:delta}
\end{equation}
In one of our previous
publications \cite{sol3} we  tried to connect the appearance of
surface terms in Poisson brackets and the standard manipulations with the
$\delta$-function. The ansatz used there for
the above simplest example coincided with (\ref{eq:delta}) up to the sign.
The reason for the
difference laid in the other choice made there instead of Rule 4.2 of I.
That ansatz led us to the standard Poisson brackets which were not
appropriate for boundary problems.

Operators satisfying relation
\[
\hat I^{\ast}=-\hat I
\]
will be called {\it skew-adjoint} ones. With the help of them it is possible
to express 2-forms (and also 2-vectors to be defined below) in the canonical
form \[ \alpha=\frac{1}{2}\sum\limits_{A,B}\int\delta\phi_A\wedge\hat I_{AB}
\delta\phi_B.
\]
It is clear that we can consider representations of functional forms
as decompositions over the basis derived as a result of the tensor product
of $\delta\phi_A$, with the totally antisymmetric multilinear operators
\[
\hat\alpha=\sum\theta^{(J)}\alpha^{\langle J\rangle}_{A_1K_1,\dots,A_mK_m}
\biggl( D_{K_1}\cdot,\dots,D_{K_m}\cdot\biggr)
\]
as coefficients of these decompositions.

\section{Multi-vectors, mixed tensors and Schouten-Nijenhuis bracket}
Let us introduce dual basis to $\vert\delta\phi_A\rangle$ by relation
\begin{equation}
\left\langle\frac{\delta}{\delta\phi_B(y)},\delta\phi_A(x)\right\rangle
=\delta_{AB}\delta(x,y) \label{eq:dual}
\end{equation}
and construct by means of the tensor product a basis
\[
\frac{\delta}{\delta\phi_{B_1}(y)}\otimes\frac{\delta}{\delta\phi_{B_2}(y)}
\otimes\dots\otimes\frac{\delta}{\delta\phi_{B_m}(y)}.
\]
Then by using totally antisymmetric multilinear operators described in
previous Section we can define {\it functional $m$-vectors} (or {\it
multi-vectors})
\[
\psi=\frac{1}{m!}\sum\int\theta^{(J)} \psi^{\langle
J\rangle}_{B_1L_1,\dots,B_mL_m}D_{L_1}
\frac{\delta}{\delta\phi_{B_1}}\wedge\dots\wedge D_{L_m}
\frac{\delta}{\delta\phi_{B_m}}.
\]
Here a natural question on the relation between evolutionary
vector fields and 1-vectors arises.
Evidently, evolutionary vector fields lose their form when being
integrated by parts whereas 1-vectors conserve it.
Let us make a partial integration in the expression
of a general evolutionary vector field
\[
\xi=\sum\int\theta^{(J)}
D_K\xi^{\langle J\rangle}_A\frac{\partial}{\partial\phi_A^{(K)}}
\]
by removing $D_K$ from $\xi^{\langle J\rangle}_A$, then we get
\[
\xi=\sum
\int\xi^{\langle J\rangle}_A\theta^{(J+L)}
(-1)^{|K|}{K \choose L}D_{K-L}\frac{\partial}{\partial\phi_A^{(K)}}.
\]
It is easy to see that by using Rule 5.4 from I in the backward direction
we can write
\[
\xi=\sum\int\bigl[ \theta^{(J)}\xi_A^{\langle J\rangle}\bigr]\bigl[
\theta^{(L)}(-1)^{|L|}E^L_A\bigr]=\sum\int\theta^{(J)}\xi^{\langle
J\rangle}_A
\frac{\delta}{\delta\phi_A},
\]
where the higher Eulerian operators and full variational derivative
(Definition 5.1 of I)
are consequently used. Therefore, we have proved a following Statement.

{\bf Statement 5.1}
{\it There is a one-to-one correspondence between evolutionary vector
fields and functional 1-vectors.
The  coefficients of 1-vector in the canonical form $\xi_A^{\langle J
\rangle}$ are equal to the characteristic
of the evolutionary vector field.}

It is not difficult to show that we can  deduce the pairing (interior
product) of 1-forms and 1-vectors and this pairing preserves this
identification.  Really, the definition of the dual basis (\ref{eq:dual}) and
Rules 4.2, 5.4 of I permit us to derive that
\[
\alpha(\xi)=\xi\inprod\alpha=\sum\int \int\theta^{(I)}(x)
\theta^{(J)}(y)\alpha^{\langle I\rangle}_{AK}(x)\xi^{\langle J\rangle
}_{BL}(y)
\left\langle D_L\frac{\delta}{\delta\phi_B(y)},D_K\delta\phi_A(x)\right
\rangle=
\]
\[
=\sum\int\theta^{(I+J)}D_L\alpha^{\langle I\rangle}_{AK}
D_K\xi^{\langle J\rangle}_{AL}=
\sum\int\theta^{(I+J)}{\rm Tr}(\alpha^{\langle I\rangle}
\xi^{\langle J\rangle}),
\]
and when 1-vector is in the canonical form (only $L=0$ term is nonzero)
this result coincides with Eq.(\ref{eq:pairing}).

This formula for the pairing will be exploited below also for interior
product of 1-vectors and $m$-forms or 1-forms and $m$-vectors.
Its importance comes from the fact that it is invariant under the formal
partial integration both in forms and in vectors, i.e., \[ {\rm
Div}(\alpha)(\xi)={\rm Div}(\alpha(\xi))=\alpha({\rm Div}(\xi)).  \]
Evidently, it is the trace construction for convolution of differential
operators (as coefficients of tensor objects in the proposed basis)
that guarantees this invariance.

The interior product of 1-vector onto $m$-form and, analogously,
of 1-form onto $m$-vector is defined as
\[
\xi\inprod\alpha=\frac{1}{m!}\sum (-1)^{(i+1)}\int
\theta^{(I+J)}D_{K_i}\xi^{\langle I\rangle}_{A_iL}D_L
\biggl(\alpha^{\langle J\rangle}_{A_1K_1,\dots,
A_mK_m}\delta\phi_{A_1}^{(K_1)}\wedge\dots
\]
\[
\dots\wedge\delta\phi_{A_{i-1}}^{(K_{i-1})}\wedge
\delta\phi_{A_{i+1}}^{(K_{i+1})}\wedge\dots
\wedge\delta\phi_{A_m}^{(K_m)}\biggr).
\]
Then we also can define the value of $m$-form on $m$ 1-vectors (or,
analogously, $m$-vector on $m$ 1-forms)
\[
\alpha (\xi_1,\dots,\xi_m)=
\xi_m\inprod\dots\xi_1\inprod\alpha=
\sum\int\theta^{(J+I_1+\dots+I_m)}{\rm Tr}
\biggl( \alpha^{\langle J\rangle}
\xi_1^{\langle I_1\rangle}\cdots\xi_m^{\langle I_m\rangle}\biggr),
\]
where each entry of multilinear operator $\alpha$ acts only to the one
$\xi$, whereas each derivation of the operator $\xi$ acts on the
product of $\alpha$ and all the rest of $\xi$'s.

It is possible to define the {\it differential of $m$-vector}
\[
\d\psi =\frac{1}{m!}\sum\int\theta^{(J)}
\frac{\partial\psi^{\langle J\rangle}_
{A_1K_1,\dots,A_mK_m}}{\partial\phi_B^{(L)}}\delta\phi_B^{(L)}D_{K_1}
\frac{\delta}{\delta\phi_{A_1}}\wedge\dots\wedge
D_{K_m}\frac{\delta}{\delta\phi_{A_m}},
\]
as an example of a mixed ${m \choose 1}$ object. Evidently, ${\d}^2\psi=0$.

With the help of the previous constructions we can define the
{\it Schouten-Nijenhuis bracket}
\[
\bigl[ \xi,\eta\bigr]_{SN} =\d\xi\inprod\eta +
(-1)^{pq}\d\eta\inprod\xi
\]
for two multi-vectors of orders
$p$ and $q$. The result of this operation is $p+q-1$-vector and it is
analogous to the Schouten-Nijenhuis bracket in tensor analysis \cite{nij}.
Its use in the formal variational calculus is described in
Refs.\cite{olv},\cite{dorf}.  However, in cited references this bracket is
usually defined for operators.  We can recommend Ref.\cite{olv2} as an
interesting source for the treatment of the Schouten-Nijenhuis bracket of
multi-vectors.  Our construction of this bracket guarantees a compatibility
with the equivalence modulo divergences
\[ \bigl[ {\rm
Div}(\xi),\eta\bigr]_{SN} ={\rm Div}\bigl[ \xi,\eta\bigr]_{SN}= \bigl[ \xi,
{\rm Div}(\eta)\bigr]_{SN}.
\]
\medskip
{\bf Statement 5.2} {\it The
Schouten-Nijenhuis bracket of functional 1-vectors up to a sign coincides
with the commutator of the corresponding evolutionary vector fields.}
\medskip

{\it Proof.}
Let us take the two 1-vectors in canonical form without loss of generality
\[
\xi=\sum\int\theta^{(J)}\xi^{\langle J\rangle}_A
\frac{\delta}{\delta\phi_A},\qquad
\eta=\sum\int\theta^{(K)}\eta^{\langle K\rangle}_B\frac{\delta}{\delta\phi_B}
\]
and compute
\[
\bigl[ \xi,\eta\bigr]_{SN}=\d\xi\inprod\eta
-\d\eta\inprod\xi.
\]
We have
\[
\d\xi=\sum\int\theta^{(J)}{\xi^{\langle
J\rangle}_A}'(\delta\phi)\frac{\delta}{\delta
\phi_A}=\sum\int\theta^{(J)}\frac{\partial\xi^{\langle
J\rangle}_A}{\partial\phi^{(L)}
_C}\delta\phi_C^{(L)}\frac{\delta}{\delta\phi_A},
\]
and
\[
\d\xi\inprod\eta=-\sum\int\theta^{(J+K)}\frac{\partial
\xi_A^{\langle J\rangle}}
{\partial\phi_B^{(L)}}D_L\eta_B^{\langle K\rangle}\frac{\delta}{\delta\phi_A}.
\]
Therefore, we obtain
\[
\bigl[ \xi,\eta\bigr]_{SN}=-\sum\int\theta^{(J+K)}\biggl(
D_L\eta_B^{\langle K\rangle}\frac{\partial\xi_A^{\langle J\rangle}}
{\partial\phi_B^{(L)}}-
D_L\xi_B^{\langle K\rangle}\frac{\partial\eta_A^{\langle J\rangle}}
{\partial\phi_B^{(L)}}
\biggr)
\frac{\delta}{\delta\phi_A}=-[\xi,\eta],
\]
and the proof is completed.

\medskip
{\bf Statement 5.3} (Olver's Lemma \cite{olv})
{\it
The Schouten-Nijenhuis bracket for two bivectors can be expressed
in the form}
\begin{equation}
\bigl[ \xi,\psi\bigr]_{SN}=-\frac{1}{2}\sum\int\
\xi\wedge\hat I'(\hat K\xi)\wedge\xi
-\frac{1}{2}\sum\int\
\xi\wedge\hat K'(\hat I\xi)\wedge\xi, \label{eq:prolong}
\end{equation}
{\it where the two differential operators $\hat I$, $\hat K$ are the
coefficients of the bivectors in their canonical form.}
\medskip

{\it Proof.}
Let us consider the Schouten-Nijenhuis bracket for the two bivectors and
without loss of generality
take them in the canonical form
\[
\chi=\frac{1}{2}\sum\int\theta^{(L)}\xi_A\wedge I^{\langle L\rangle N}_{AB}
D_N\xi_B,
\]
\[
\psi=\frac{1}{2}\sum\int\theta^{(M)}\xi_C\wedge K^{\langle M\rangle P}_{CD}
D_P\xi_D,
\]
where $\xi_A={\delta}/{\delta\phi_A}$ and operators $\hat I$ , $\hat K$
are skew-adjoint. Then we have
\[
\d\chi=\frac{1}{2}\sum\int\theta^{(L)}\frac{\partial
I^{\langle L\rangle N}_{AB}}
{\partial\phi_E^{(J)}}\delta\phi_E^{(J)}\xi_A\wedge D_N\xi_B
\]
and
\[
\d\chi\inprod\psi=\frac{1}{4}\sum\int\theta^{(L+M)}
\frac{\partial I^{\langle L\rangle N}_{AB}}{\partial\phi_C^{(J)}}D_J
\biggl( K^{\langle M\rangle P}_{CD}D_P\xi_D\biggr)\wedge
\xi_A\wedge D_N\xi_B-
\]
\[
-\frac{1}{4}\sum\int\theta^{(L+M)}D_P\biggl(
\frac{\partial I^{\langle L\rangle N}_{AB}}
{\partial\phi_D^{(J)}}\xi_A\wedge D_N\xi_B\biggr)
\wedge D_J(\xi_C K^{\langle M\rangle P}_{CD}).
\]
Now let us make integration by parts in the second term
\[
\d\chi\inprod\psi=-\frac{1}{4}\sum\int\theta^{(L+M)}\xi_A\wedge
(I^{\langle L\rangle N}_{AB})'\biggl(\hat K^{\langle M\rangle}
\xi\biggr)\wedge  D_N\xi_B-
\]
\[
-\frac{1}{4}\sum\int\theta^{(L+M+Q)}(-1)^{|P|}{P\choose Q}
\frac{\partial I^{\langle L\rangle N}_{AB}}
{\partial\phi_D^{(J)}}\xi_A\wedge D_N\xi_B
\wedge D_{J+P-Q}(\xi_C K^{\langle M\rangle P}_{CD}).
\]
At last we change the order of multipliers under wedge product
in the second term,
make a replacement $M\rightarrow M-Q$ and organize the
whole expression in the form
\[
\d\chi\inprod\psi=-\frac{1}{4}\sum\int\theta^{(L+M)}\xi_A\wedge
(I^{\langle L\rangle N}_{AB})'_C\Biggl(\hat K^{\langle M\rangle}_{CD}\xi_D-
\]
\[
-(-1)^{|P|}{P\choose Q}{P-Q\choose R}
D_{P-Q-R}K^{\langle M-Q\rangle P}_{CD}D_R\xi_C \Biggr)\wedge D_N\xi_B.
\]
Having in mind the definition of adjoint operator (\ref{eq:adj}) we can
represent the final result of the calculation as follows,
\[
\bigl[ \xi,\psi\bigr]_{SN}=-\frac{1}{2}\sum\int\theta^{(L+M)}
\xi\wedge\biggl((\hat I^{\langle L\rangle})'(\hat K^{\langle M\rangle}\xi)
-(\hat K^{\langle M\rangle})'(\hat I^{\langle L\rangle}\xi)\biggr)\wedge\xi,
\]
thus supporting in this extended formulation the method, proposed
in Ref.~\cite{olv} for testing the Jacobi identity (see Section 7).

\section{Poisson brackets and Hamiltonian vector fields}
Let us call a bivector
\[
\Psi=\frac{1}{2}\sum\int\frac{\delta}{\delta\phi_A}\wedge\hat I_{AB}
\frac{\delta}{\delta\phi_B},
\]
formed with the help of the graded skew-adjoint differential operator
\[
\hat I_{AB}=\sum \theta^{(L)}I^{\langle L\rangle N}_{AB}D_N,
\]
the {\it Poisson bivector} if
\[
\bigl[ \Psi,\Psi\bigr]_{SN} =0.
\]
The operator $\hat I_{AB}$ is then called the {\it Hamiltonian operator}.
We call the value of the Poisson bivector on the differentials of
two functionals $F$ and  $G$
\[
\{ F,G\} = \Psi (\d F,\d G)=\d G\inprod\d F\inprod\Psi
\]
the {\it Poisson bracket} of these functionals.

The explicit form of the Poisson brackets can easily be obtained. It
depends on the explicit form of the functional differential,
which can be changed by partial integration. Of course, all the
possible forms are equivalent. Taking the extreme cases we get an
expression through Fr\'echet derivatives
\begin{equation}
\{ F,G \} =
\sum\int\theta^{(J)} {\rm Tr}\biggl( f'_A\hat I^{\langle J\rangle}_{AB}
g'_B \biggr) \label{eq:brack1}
\end{equation}
or through higher Eulerian operators
\begin{equation}
\{ F,G \} =
\sum\int\theta^{(J)} D_{P+Q}\biggl( E^P_A(f)\hat I^{\langle J\rangle}_{AB}
E^Q_B(g)
\biggr).\label{eq:brack2}
\end{equation}

\medskip
{\bf Theorem 6.1}
{\it The Poisson bracket defined above satisfies}
Definition 2.3 {\it of} I.
\medskip

{\it Proof.}
Equivalence of these definitions  follows from the three
facts: 1) from the previous formulas (\ref{eq:brack1}), (\ref{eq:brack2})
it is clear that $\{ F,G \}$ is a local functional,
2) antisymmetry of $\{ F,G \}$ is evident and
3) equivalence of the Jacobi identity
to the Poisson bivector property (to be proved in Section 7).
\medskip

The result of interior product of the differential of a local
functional $H$ on the Poisson bivector (up to the sign)
will be called the {\it  Hamiltonian
vector field} (or the {\it Hamiltonian 1-vector})
\[
\hat I\d H=-\d H\inprod\Psi
\]
corresponding to the Hamiltonian $H$.
Evidently, the standard relations take place
\[
\{ F,H\} = \d F(\hat I \d H)=(\hat I \d H)F.
\]

\medskip
{\bf Theorem 6.2}
{\it The Hamiltonian vector field corresponding to the Poisson bracket
of the functionals $F$ and $H$ coincides up to the sign with
commutator of the Hamiltonian vector fields corresponding to these
functionals.}
\medskip

{\it Proof.}
Consider a value of the commutator of Hamiltonian vector fields
$\hat I\d F$ and $\hat I\d H$ on the arbitrary functional $G$
\[
[\hat I\d F, \hat I\d H]G=\hat I\d F(\hat I\d H(G))-\hat I\d H(\hat I\d F(G))
=
\]
\[
=\hat I\d F(\{G,H\})-\hat I\d H(\{G,F\})=\{\{G,H\},F\}-\{\{G,F\},H\}=
\]
\[
=-\{G,\{F,H\}\}=-\hat I\d\{F,H\}(G),
\]
where we have used the Jacobi identity and antisymmetry of Poisson bracket.
Due to the arbitrariness of $G$ the proof is completed.

\medskip
{\it Example 6.3}

Let us consider a first structure
\[
\{ u(x),u(y)\} =\frac{1}{2}(D_x-D_y)\delta(x,y)
\]
of the Korteweg-de Vries equation (Example 7.6 of Ref.~\cite{olv})
\[
u_t=u_{xxx} +uu_x.
\]
Construct the adjoint graded operator
to $\theta D$ according to Eq.(\ref{eq:adj})
\[
(\theta D)^{\ast}=-\theta D -D\theta
\]
and the skew-adjoint operator is
\[
\hat I=\frac{1}{2}\biggl(\theta D- (\theta D)^{\ast}\biggr)=\theta D +
\frac{1}{2}D\theta.
\]
The Poisson bivector has a form
\[
\Psi=\frac{1}{2}\int\theta\biggl( \frac{\delta}{\delta u}\wedge
D\frac{\delta}{\delta u}\biggr).
\]
The differential of a local functional $H$ (for simplicity
it is written in canonical
\[
H=\int\theta h
\]
form) is equal to
\[
\d H=\int\theta h'(\delta u)=\sum_{k=0}^{\infty}
\int\theta^{(k)}(-1)^kE^k(h)\delta u,
\]
where the Fr\'echet derivative or higher Eulerian operators can be used.
Therefore, the Hamiltonian vector field generated by $H$ is
\[
\hat I\d H=-\d H\inprod\Psi=
-\frac{1}{2}\int\theta\biggl[ h'\bigl( D\frac{\delta}{\delta u}
\bigr) - Dh'\bigl( \frac{\delta}{\delta u}\bigr)\biggr] ,
\]
or
\[
-\frac{1}{2}\int\theta^{(k)}(-1)^k\biggl[ E^k(h)D-
DE^k(h)\biggr]\frac{\delta}{\delta u},
\]
or also
\[
-\frac{1}{2}\int\theta^{(k)}(-1)^kD_i\biggl[ E^k(h)D-
DE^k(h)\biggr]\frac{\partial}{\partial u^{(i)}}.
\]
The value of this vector field on another functional $F$ coincides with
the Poisson bracket
\[
-\d F\inprod\d H\inprod\Psi=\{ F,H\}=
\frac{1}{2}\sum\int\theta^{(k+l)}(-1)^{k+l}
\biggl( E^k(f)DE^l(h)-E^k(h)DE^l(f)\biggr).
\]

\section{Proof of Jacobi identity}
In this section we will prove that the Jacobi identity for the Poisson bracket
is fulfilled if and only if the Schouten--Nijenhuis bracket of the
corresponding Poisson bivector with itself is equal to
zero. This should complete the proof of Theorem 6.1.

Let us use one of the possible forms of the Poisson brackets given in
Appendix of I
\[
\{ F,G\} =\frac{1}{2}\sum\int\theta^{(J)}{\rm Tr}\biggl( f'(\hat I^{\langle
J\rangle}g')
-g'(\hat I^{\langle J\rangle}f')\biggr),
\]
where the differential operator $\hat I$ is not supposed to be skew-adjoint
for the easier comparison of this proof with that given in I.
We remind the reader that in less condensed notations
\[
{\rm Tr}\biggl( f'(\hat Ig')\biggr)=\sum{J\choose M}{K\choose L}
D_L\frac{\partial f}{\partial\phi_A^{(J)}}D_{J+K-L-M}I^N_{AB}
D_{N+M}\frac{\partial g}{\partial\phi_B^{(K)}}
\]
(in Appendix of I the indices $M$ and $L$ in the binomial coefficients
of the same formula  are unfortunately given in the opposite order).

We will estimate the bracket
\[
\{\{ F,G\} ,H\} =\frac{1}{2}\sum\int\theta^{(J)}{\rm Tr}\biggl[
{\{ f,g\} }'(\hat I^{\langle J\rangle} h')-h'(\hat I^{\langle J\rangle}
{\{ f,g\}}')\biggr],
\]
where $\{ f,g\}$ denotes the integrand of $\{ F,G\}$. Since Fr\'echet
derivative is a derivation we have
\[
{\{ f,g\} }'=\frac{1}{2}\sum\theta^{(K)}{\rm Tr}\biggl(
f''(\hat I^{\langle K\rangle}g',\cdot)
+f'\hat I'^{\langle K\rangle}(\cdot)g'+
g''(f'\hat I^{\langle K\rangle},\cdot)-(f\leftrightarrow g)
\biggr)
\]
and
\[
{\rm Tr}\biggl[ {\{ f,g\}}'\hat Ih'\biggr]=\frac{1}{2}\biggl[
f''(\hat Ig',\hat Ih')+f'\hat I'(\hat Ih')g'+g''(f'\hat I,\hat Ih')-
(f\leftrightarrow g)\biggr].
\]
Let us explain that $f''$ denotes the second Fr\'echet derivative, i.e.,
the symmetric bilinear operator arising in calculation of the second
variation of the local functional $F$ (in the canonical form):
\[
f''(\xi,\eta)=\sum\limits_{A,B}\sum\limits_{J,K}\frac{\partial^2f}
{\partial\phi_A^{(J)}\partial\phi_B^{(K)}}D_J\xi_AD_K\eta_B.
\]
When we  put into entries of $f''$ operators under the trace sign
it should be understood that these operators act on everything except
their own coefficients, for example,
\[
{\rm Tr}\biggl(f''(\hat Ig',\hat Ih')\biggr)
=\sum{L\choose P}{L-P\choose Q}{M\choose S}
{M-S\choose T}\times
\]
\[
\times D_{L+M-P-Q-S-T}\frac{\partial^2f}{\partial\phi_A^{(J)}
\partial\phi_B^{(K)}}
D_{J+T}\biggl( D_P\hat I_{AC}\frac{\partial g}
{\partial\phi_C^{(L)}}\biggr) D_{K+Q}\biggl( D_S\hat I_{BD}\frac{\partial h}
{\partial\phi_D^{(M)}}\biggr)
\]
and the expression remains symmetric under permutation of its entries
\[
{\rm Tr}\biggl( f''(\hat Ig',\hat Ih')\biggr)=
{\rm Tr}\biggl( f''(\hat Ih',\hat Ig')\biggr) .
\]
When the operator $\hat I$  stands to the right from the operator
of Fr\'echet derivative $f'$ as in expression
\[
{\rm Tr}\biggl( g''(\hat Ih',f'\hat I)\biggr) ,
\]
it acts on everything except $f'$. At last, for Fr\'echet derivative
of the operator we have
\[
\hat I'(\hat Ih')=\sum\frac{\partial I^K_{AB}}{\partial\phi_C^{(J)}}
D_J\biggl( I^L_{CD}D_L\frac{\partial h}{\partial\phi_D^{(M)}}D_M\biggr) D_K.
\]
Making similar calculations we get
\[
{\rm Tr}\biggl[ h'\hat I{\{ f,g\}}'\biggr]=
\frac{1}{2}{\rm Tr}\biggl( f''(h'\hat I,\hat Ig')+f'\hat I'(h'\hat I)g'+
g''(f'\hat I,h'\hat I)-(f\leftrightarrow g)\biggr)
\]
and therefore
\[
\{\{ F,G\} ,H\} =\frac{1}{4}\sum\int\theta^{(J+K)}{\rm Tr}\biggl(
f''(\hat I^{\langle J\rangle}g',\hat I^{\langle K\rangle}h')-
f''(h'\hat I^{\langle J\rangle},\hat I^{\langle K\rangle}g')-
\]
\[
-f''(\hat I^{\langle J\rangle}h',g'\hat I^{\langle K\rangle})+
f''(g'\hat I^{\langle J\rangle},h'\hat I^{\langle K\rangle})+
f'\hat I'^{\langle J\rangle}(\hat I^{\langle K\rangle}h'-h'
\hat I^{\langle K\rangle})g' -(f\leftrightarrow g)
\biggr) .
\]
Just the first four terms, apart from the fifth containing Fr\'echet
derivative of the
operator $\hat I$, were present in our proof for nonultralocal case
given in I (only terms with zero grading were allowed for $\hat I$ there).
After cyclic permutation of $F$, $G$, $H$ all terms with the symmetric
operator of the second Fr\'echet derivative are mutually cancelled and
\[
\{\{ F,G\} ,H\} + {\rm c.p.}=\frac{1}{4}\int\theta^{(J+K)}{\rm Tr}\biggl(
f'\hat I'^{\langle J\rangle}(\hat I^{\langle K\rangle}h'-
h'\hat I^{\langle K\rangle })g'-
\]
\[
-g'\hat I'^{\langle J\rangle}
(\hat I^{\langle K\rangle}h'-h'\hat I^{\langle K\rangle})f'+ c.p.
\biggr),
\]
where cyclic permutations of $F$, $G$, and $H$ are abbreviated to $c.p.$.
When operator $\hat I$ is given in explicitly skew-adjoint form all
the four terms are equal. Taking into account Olver's Lemma
(\ref{eq:prolong})
we get
\[
\{\{ F,G\} ,H\} +  c.p. =-\bigl[\hat I,\hat
I\bigr]_{SN}(\d F,\d G,\d H),
\]
so finishing the proof.

\section{Examples of nonultralocal operators}
The second structure of the Korteweg-de Vries equation may serve as a
counter-example to hypothesis \cite{coll}
that all operators which are Hamiltonian
with respect to the standard Poisson brackets should also be Hamiltonian in
the new brackets.

\medskip
{\it Example 8.1}

Let us start with the standard expression (Example 7.6 of Ref.~\cite{olv})
\[
\{ u(x),u(y)\} =\biggl(\frac{d^3}{dx^3}+\frac{2}{3}u\frac{d}{dx}+
\frac{1}{3}\frac{du}{dx})\delta(x,y)
\]
and construct the adjoint operator to
\[
\hat K=\theta(D_3+\frac{2}{3}uD+\frac{1}{3}Du),
\]
which is
\[
\hat K^{\ast}=-\theta(D_3+\frac{2}{3}uD+\frac{1}{3}Du)-D\theta(3D_2+
\frac{2}{3}u)-3D_2\theta D-D_3\theta.
\]
Then the skew-adjoint operator
\[
\hat I=\frac{1}{2}(\hat K-\hat K^{\ast})=
\theta(D_3+\frac{2}{3}uD+\frac{1}{3}Du)+D\theta(\frac{3}{2}D_2+\frac{1}{3}u)
+\frac{3}{2}D_2\theta D+\frac{1}{2}D_3\theta
\]
can be used for forming the bivector
\[
\Psi=\frac{1}{2}\int\xi\wedge\hat I\xi,
\]
where  ${\delta}/{\delta u}=\xi$. This bivector has a form
\[
\Psi=\frac{1}{2}\int\biggl(\theta\xi\wedge D_3\xi+\frac{3}{2}D\theta\xi
\wedge D_2\xi+(\frac{3}{2}D_2\theta+\frac{2}{3}\theta u)\xi\wedge D\xi\biggr).
\]
Then evaluating the Schouten-Nijenhuis bracket for the bivector with
the help of Statement 5.3
\[
\bigl[\Psi,\Psi\bigr]_{SN}=\int (\frac{2}{3}\theta\xi\wedge D_3\xi
\wedge D\xi +
D\theta\xi\wedge D_2\xi\wedge D\xi)
\]
and  integrating the first term by parts we get
\[
\bigl[\Psi,\Psi\bigr]_{SN}=\frac{1}{3}\int\theta D\bigl(\xi\wedge D\xi
\wedge D_2\xi\bigr).
\]
Therefore, instead of the Jacobi identity we  have
\[
\{\{ F,G\},H\}+c.p.=-\frac{1}{3}\sum\limits_{i,j,k=0}^{\infty}
\int\limits_{\Omega}D_{i+j+k+1}\biggl( E^i(f)DE^j(g)D_2E^k(h)+c.p.\biggr)dx.
\]
So, the second structure of KdV equation can be Hamiltonian only
under special boundary conditions.

\medskip
{\it Example 8.2}

Now let consider another example which is also nonultralocal, but the operator
remains
to be Hamiltonian in the new brackets independently of boundary conditions.
The Euler equations for the flow of ideal  fluid
can be written \cite{olv} in Hamiltonian form as follows
(Example 7.10 of Ref.~\cite{olv})
\[
\frac{\partial{\bf\omega}}{\partial t}={\cal D}\frac{\delta H}
{\delta{\bf\omega}},
\]
where
\[
H=\int\frac{1}{2}\vert {\bf u}\vert^2d^2x,\qquad {\bf\omega}={\bf\nabla}
\times{\bf u}.
\]
Let us limit our consideration by the 2-dimensional case when $\bf\omega$
has only one component $\omega$ and
\[
{\cal D}={\bf\omega}_xD_y-{\bf\omega}_yD_x,
\]
where $\omega_i=D_i\omega$, $i=(x,y)$.
We can construct the skew-adjoint operator
\[
\hat I=\frac{1}{2}\biggl( \theta{\cal D}-(\theta{\cal D})^{\ast}\biggr)=
\theta(\omega_xD_y-\omega_yD_x)
+\frac{1}{2}(D_y\theta\omega_x-D_x\theta\omega
_y),
\]
and then the bivector
\[
\Psi=\frac{1}{2}\int\xi\wedge\hat I\xi=
\frac{1}{2}\int\theta(\omega_x\xi\wedge
\xi_y-\omega_y\xi\wedge\xi_x),
\]
where $\xi={\delta}/{\delta\omega}$.
Statement 5.3 gives us
\[
\bigl[\Psi,\Psi\bigr]_{SN}=
\int\Biggl(\theta\biggl[
\omega_x(\xi\wedge\xi_{xy}\wedge\xi_y-\xi\wedge\xi_{yy}
\wedge\xi_x)+\omega_y(\xi\wedge\xi_{xy}\wedge\xi_x-\xi\wedge\xi_{xx}\wedge
\xi_y)\biggr]+
\]
\[
+\biggl[ D_y\theta\omega_x-
D_x\theta\omega_y\biggr]\xi\wedge\xi_x\wedge\xi_y\Biggr)
\]
and after integration by parts the expression can be reduced to zero.

\section{Conclusion}
We have shown that there is an extension of the standard formal
variational calculus which incorporates the real divergences without
any specification of the boundary conditions on the boundary of a
finite domain. It would be important to understand relations of this
formalism to the constructions of the variational bicomplex \cite{and}.
It seems also rather interesting to study if some physically relevant
algebras can be realized with the help of the new Poisson brackets
as algebras of local functionals. It is not clear for us now whether
the Hamiltonian equations generated by the new brackets can be solved in
some space of functions and what kind of space could be used for this
purpose.

\vspace{12pt}

{\large\bf Acknowledgements}

It is a pleasure to thank S.N.Storchak for discussions and A.V.Razumov
for answering numerous questions.

This work has been started during a visit of the author to the International
Centre for Theoretical Physics in Trieste, partial support from ICTP is
gratefully acknowledged.

\hfill

\end{document}